\documentclass[aps,pre,twocolumn,showpacs,groupedaddress,amssymb]{revtex4}
\usepackage{graphicx,float,subfig,amsmath,booktabs}% Include figure files

\begin{document}

\title{Partial differential equations versus cellular automata for modelling combat}

\author{Therese Keane}
\altaffiliation[Permanent address: ]{Defence Science and Technology Organisation, Eveleigh, Sydney NSW 2015, Australia.}
\email[]{Therese.Keane@dsto.defence.gov.au}
\affiliation{Department of Mathematics and Statistics, University of New South Wales, Sydney NSW 2052, Australia.}

\date{\today}

\begin{abstract}

We reproduce apparently complex cellular automaton behaviour with simple partial differential equations as developed in \cite{tk09}.  Our PDE model easily explains behaviour observed in selected scenarios of the cellular automaton wargame ISAAC without resorting to anthropomorphisation of autonomous 'agents'.  The insinuation that agents have a reasoning and planning ability is replaced with a deterministic numerical approximation which encapsulates basic motivational factors and demonstrates a variety of spatial behaviours approximating the mean behaviour of the ISAAC scenarios.  All scenarios presented here highlight the dangers associated with attributing \emph{intelligent} reasoning to behaviour shown, when this can be explained quite simply through the effects of the terms in our equations.  A continuum of forces is able to behave in a manner similar to a collection of individual autonomous agents, and shows decentralised self-organisation and adaptation of tactics to suit a variety of combat situations.

\end{abstract}

%\pacs{87.18.-h,02.60.-x,87.75.Kd}
%\keywords{}

\maketitle

\section{Introduction\label{sec:Introduction}}
Lanchester modelled force dynamics as a set of ordinary differential equations that have greatly influenced military decision making for many years and permeate military thought and analysis to this day.  However, the absence of any spatial representation has been a constant source of criticism for using any form of the Lanchester Equations.  Extensions have been made by Protopopescue \emph{et al.} who treated the spatial aspect using a simple advection-diffusion representation as is typical in predator-prey type modelling \cite{p87,p89,p90}.  Despite using specific initial density profiles and short simulation times to prevent excessive and unrealistic diffusion of the forces, this approach was able to represent a suite of realistic behaviour - frontal attack, envelopment, infiltration and turning manoeuvres.  However, replacing the constant velocity with a spatially dependent velocity field resulted in unacceptable numerical losses, restricting the velocity field to a temporally dependent one.

  In \cite{tk09} a set of partial differential equations for a more spatially realistic approach to combat modelling was proposed.  Through the implementation of biological aggregation models, a series of basic scenarios demonstrating cohesive soldier movement throughout the domain regardless of losses incurred through firing effects were explored.  It was also shown that a scenario can be viewed as a series of sub-battles characterised by periods of constant density loss which was of the same form as in \cite{p87,p89,p90}.  This was the intended method of implementation of the LEs as suggested by Lanchester himself.  These PDEs were developed with the intention of providing an alternate method for investigating the dynamics seen in cellular automaton-based wargames.
  
	These agent-based or cellular-automaton (CA) models have received much attention in many disciplines, particularly in defence related research.  Models such as Einstein \cite{andy03}, ISAAC \cite{andy97}, Map Aware Non-uniform Automata (MANA) \cite{l02} demonstrate a range of behaviour which appears to hint at some form of underlying structure.  Each individual troop is modelled via a rule set relating to quantifiable capabilities such as fire-power, communications and also intangibles such as morale or desire to remain close to friendly forces.  These rules encode the nonlinearities necessary for a more realistic description of warfare.  These nonlinearities need to be understood in order to develop specialised tactics based on current capability, or enhance the procurement of future capability.  Ilachinski \cite{andy96}, who has been instrumental in the development of ISAAC, stresses the need for research into nonlinear continuous dynamics, exploitation of analogous biological models and phase-space reconstruction techniques.  Lauren compares MANA simulation results with fluid dynamic concepts or transition between laminar and turbulent states and maintenance of force profiles to viscosity \cite{l99}.
	
	Continuous models can be more transparent in terms of how parameter changes affect outcomes and are thus more understandable.  In this paper, the PDE model is expanded to reflect the underlying assumptions of the ISAAC model to provide an alternate analysis of the demonstrated behaviour.  By effectively removing the inherent stochasticity of ISAAC, we can determine the mean behaviour of a scenario and more easily understand the effects of parameter changes.  Introducing spatial asymmetry through the initial density profile position forms a type of controlled stochasticity approximation, providing an explanation of the observed ISAAC dynamics.
\subsection{PDE model}
As developed in \cite{tk09}, we begin with the following integro-differential equations in two dimensions.  Each equation is separated into three parts, $\mathbf{f_{diff}}$ containing the diffusion terms, $\mathbf{f_{vel}}$ containing the velocity terms and $\mathbf{f_{react}}$ containing the interaction terms.
\begin{align}\label{eq:NewEquationForU}
	\frac{\partial u}{\partial t} &= \underbrace{\nabla \cdot(\mathbf{D_u}(u)\nabla u)}_{\mathbf{f_{diff}}} + \underbrace{u\left(k_{u} \ast v\right) + d_uv}_{\mathbf{f_{react}}} +\nonumber \\
	& \qquad \underbrace{\nabla \cdot\{u(\mathbf{C_u}u + A_{a}(K_{a} \ast u) - A_{r}u(K_{r} \ast u))\}}_{\mathbf{f_{vel}}}
\end{align}
\begin{align}\label{eq:NewEquationForV}
	\frac{\partial v}{\partial t} &= \underbrace{\nabla \cdot(\mathbf{D_v}(v)\nabla v)}_{\mathbf{f_{diff}}}+ \underbrace{v\left(k_{v} \ast u\right) + d_vu}_{\mathbf{f_{react}}} + \nonumber \\
	& \qquad \underbrace{\nabla \cdot\{v(\mathbf{C_v}v + A_{a}(K_{a} \ast v) - A_{r}v(K_{r} \ast v))\}}_{\mathbf{f_{vel}}}
\end{align}
	Attraction and repulsion operate over different spatial domains, $r_a$ and $r_r$ respectively with $r_a > r_r$, and a spatially dependent vector field $C$ is used.  The form of the kernel in $\mathbf{f_{react}}$ is $k(x,y) = \beta e^{-\nu \sqrt{\left|((x-X)^2 + (y-Y)^2)-r_{op}\right|} } : \mathbb{R}^2 \rightarrow \mathbb{R}^+$.
	
	Following the work of Boswell \cite{b03} and Hundsdorfer \cite{hun95}, the Method of Lines scheme is implemented where the spatial derivatives are discretised, giving a large system of Ordinary Differential Equations.  Flux limiters are employed at this stage to ensure positivity and conservation of mass.  An explicit Runge-Kutta method is then used for the time integration, again with constraints in place to ensure positivity.
\section{ISAAC}
In this section the cellular automata model ISAAC and a selection of scenarios its authors have published to demonstrate its range of capabilities are described.  We then show how our simulation with slight modifications can reproduce similar behaviours.  Thus we show that a PDE model can perform the combat modelling tasks of a state-of-the-art CA model while being easier to understand and analyse.  This can have great advantages when seemingly unexplainable or novel behaviour is seen in CA modelling results that may otherwise be attributed to a form of intelligence.

	ISAAC, Irreducible Semi-Autonomous Adaptive Combat, is a multi-agent based simulation in the style of a cellular automaton.  A CA is essentially a lattice where information located on the nodes propagates from its position at each time step based on a set of defined rules.  Ilachinski asserts ISAAC differs from a traditional CA in that agents rather than information move throughout the lattice and that rule sets may adapt over time.  In the comparison ISAAC scenarios, however, rule sets remain fixed.
	
	In ISAAC, each agent represents a simplified soldier.  Soldiers are collectively grouped into ``Forces'' and act according to a user defined rule set for that force.  That is, all agents of a Force are \emph{homogeneous} with respect to their rule set definition.  In order to create a truly \emph{heterogeneous} battlefield of agents, each agent must be assigned to a Force consisting only of that agent.  As this could result in the creation of a large numbers of Forces and is not generally practical, two opposing Forces are usually used.  Additional to the rule sets, interaction in the form of weapons effects determines agent attrition.  Agents are defined to exist in three health states, alive, injured and killed and can obviously only exist in one of these states.  For each time step the position for each agent is determined according to the application of the previous time step data to the rule set.  Following this positional update, attrition is calculated based on the defined hit/kill probabilities.  Attrition is  applied to a specified maximum number of enemy agents (from zero to all agents) within a specified range.
	
	Movement of ISAAC agents (ISAACAs) is determined by minimisation of a penalty function calculated at each time step.  Proximity of other ISAACAs, both friendly and enemy, and proximity of each force's goals are defined as the six \emph{personality weights} used to calculate the penalty at each location within the ISAACA's movement range $r_M$.  Note that the movement range and all other ranges used in ISAAC, is a square, not a circle, of radius $r_M$.
	
	Own force goals $w_5$ are usually located at the starting side or corner of the domain/battlespace of their respective force, with the opposition force goals $w_6$ located at the opposite side.  As the proximity of an agent to the opposition goal increases, so the effect of the $\mathbf{w_6}$ term increases, increasing the speed at which it approaches the goal.  Usually the value of $\mathbf{w_5}$ is set to zero, thus the own goal has no effect of an agent.
	\begin{figure}[H]
		\centering
			\includegraphics[width=6cm]{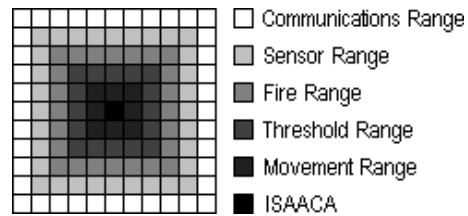}
		\caption{Ranges used in ISAAC}
		\label{fig:ISAACRanges}
	\end{figure}
Personality weights:
\begin{description}
	 \item[$\mathbf{w_1}$] alive friendly
	 \item[$\mathbf{w_2}$] alive enemy
	 \item[$\mathbf{w_3}$] injured friendly
	 \item[$\mathbf{w_4}$] injured friendly
	 \item[$\mathbf{w_5}$] friendly flag
	 \item[$\mathbf{w_6}$] enemy flag
\end{description}
The form of the penalty function is:
\begin{widetext}
\begin{align}\label{eq:ISAACPenaltyFunction}
	Z(x,y) &= \frac{w_1}{\sqrt{2}r_{f}N_{f_{alive}}}\sum_{f_{alive};i}d[i;(x,y)] + \frac{w_2}{\sqrt{2}r_{e}N_{e_{alive}}}\sum_{e_{alive};i}d[i;(x,y)] +  \frac{w_3}{\sqrt{2}r_{f}N_{f_{inj}}}\sum_{f_{inj};i}d[i;(x,y)] +\nonumber \\
	& \qquad \frac{w_4}{\sqrt{2}r_{e}N_{e_{inj}}}\sum_{e_{inj};i}d[i;(x,y)] + w_5\frac{d_{new}[flag_{f}:(x,y)]}{d_{old}[flag_{f}:(x,y)]} + w_6\frac{d_{new}[flag_{e}:(x,y)]}{d_{old}[flag_{e}:(x,y)]}
\end{align}
\end{widetext}
where $w_i, i=1,6$ are personality weights as mentioned.  The scaling factors $\sqrt{2}r_{f}$ and $\sqrt{2}r_{e}$, number of ISAACAs $N_i$ within sensor range and distances $d[i;(x,y)]$ of those $N_i$ ISAACAs form a discrete convolution.  $d_{new}$ and $d_{old}$ are the distances to the flags in question from each potential new position and from the original position respectively.

	At each time step this penalty is calculated for all potential moves an agent may make, including the penalty for remaining in the original position.  When multiple new positions of equal minimum penalty occur, the new position is chosen randomly from this set.  Interestingly, Ilachinski refers to these six personality weights as constituting a \emph{local} rule set, yet the influence of these weights, especially $w_5$ and $w_6$, can span the entire domain.
\subsection{Meta-Personality}
In addition to this penalty function, there are six additional rules that may be implemented, Advance Constraint, Cluster Constraint, Combat Constraint, Minimum distance to friendly ISAACAs, Minimum distance to enemy ISAACAs, and Minimum distance to own flag.  They are collectively termed a 'meta-personality' and modify the calculation of the penalty function.  These are effectively variations on some of the personality weights of Eq \ref{eq:ISAACPenaltyFunction} using user defined threshold and constraint ranges.

	These meta-personalities will mean our PDEs as given in (\ref{eq:NewEquationForU}) and (\ref{eq:NewEquationForV}) will require modification in order to take these behaviours into account if we wish to make a direct comparison between our model and ISAAC.  As these traits are not diffusion or reaction based, our $\mathbf{f_{vel}}$ term will be the section modified.
\subsubsection{Advance Constraint}
This constraint consists of specifying a threshold number of friendly ISAACAs that must be within a given ISAACA's constraint range $r_C$ in order for that ISAACA to continue advancing toward the enemy flag.  If the number of friendly forces within the range exceed this threshold, the default weight $+w_6$ is used in Eq \ref{eq:ISAACPenaltyFunction}.  If this is not the case, $-w_6$ is used so that the overall desire switches to a movement away from the enemy goal.  Thus sufficient troop density is required in order for advancement.
\subsubsection{Cluster Constraint}
In order for the Advance constraint to be effective, a desire to be attracted toward friendly ISAACAs to form clusters is required.   Again, the threshold number of ISAACAs required to be present within the constraint range $r_C$ is defined by the user.  If this threshold is not met, the ISAACA will move in the direction of the highest density of friendly forces calculated within the sensor range.  Once this threshold is reached, an ISAACA will no longer move toward friendly ISAACAs, effectively setting the parameters $w1 = w3 = 0$.

	Note that this constraint combined with $w_1$ and $w_3$, mimics the attraction/repulsion kernels in our continuous equations.  If the cluster constraint is not activated by the user, an \emph{artificial} default lattice repulsion is present.  Due to the density limitations of the lattice, only one agent of any type may occupy a lattice site at any time step, so that troop density cannot contract to a higher density than one troop per cell.  In order to maintain an average inner density of less than one troop per cell, the cluster constraint must be used in conjunction with $w_1$ and $w_3$.  
\subsubsection{Combat Constraint\label{sec:ISAACCombatConst}}
The Combat constraint is conceptually very similar to the Advance constraint, although this constraint is concerned with defining minimum conditions for engaging in combat with enemy ISAACAs.  Two ranges are used to determine the number of friendly and enemy ISAACAs, $N_{friendly}(r_C)$ and $N_{enemy}(r_S)$ respectively.  For advancement in the direction of greatest enemy concentration within the range $(r_S)$, the threshold troop difference $\Delta_c = N_{friendly}(r_C) - N_{enemy}(r_S)$ must be exceeded.  If this threshold is not exceeded, movement is determined using $w_2 = - w_{2,default}$ and $w_4 = - w_{4,default}$, where $w_{2,default}$ and $w_{4,default}$ are the default weights for moving toward alive and injured enemy ISAACAs.  That is, when a numerical advantage is reached, the agents will advance along the gradient of highest enemy concentration, otherwise a numerical advantage has not been reached and the agents will retreat down the gradient of highest enemy concentration.

	Setting the combat threshold to a large positive number gives a very defensive force, whereas for a large negative value a force will pursue the enemy despite the numerical disadvantage.
\subsubsection{Other Meta-Personality and Agent Constraints}
Genetic algorithms, local command and other constraints may be incorporated into ISAAC and are not discussed in this research.  Published example scenarios either do not utilise these constraints or are unaffected upon their removal.
\section{Comparisons with ISAAC Scenarios}
During the prosecution of this research, one major shortcoming of using an agent based model of this nature was highlighted - stochasticity.  Those results shown in the figures contained within the scenario descriptions of \cite{andy97} are not necessarily indicative of any expected behaviour.  Although there is some facility for the collection of basic statistics such as average cluster size, spatial entropy or Red/Blue interpoint distance, there has been no undertaking to establish whether the given realisations are true representatives of their scenario parameters.  Graham and Moyeed \cite{gm02} note that Langrangian models are akin to experiments rather than any theoretical undertaking and provide a framework for establishing reliability of such results.  We do not pursue the application of such a framework here.  Rather we take the ISAAC results presented at face value while noting the frequency of realisations that match the stated behaviour, and through our continuous approach provide a type of verification in a similar vein to \cite{go05}.
	All simulations use a domain of size $100 \times 100$ grid points.
\subsection{Classic Fronts Scenario\label{sec:ISAACClassicFronts}}
\subsubsection{Description}
We begin with the input file included with the Einstein Test Release Version 1.0.0.4 Beta, Build Date 2000, for the ``einstein\_classic\_fronts'' scenario which Ilachinski likens to a \textit{clash between two viscous fluids}.  The two loosely grouped forces collide and align in relatively stable long thin fronts.  Attrition between these fronts and the inherent randomness of the movement updates causes a discrepancy in density at either the upper or lower point of the formations.  Once this occurs, the forces are able to slowly filter around each other and proceed to their respective goals.
\subsubsection{ISAAC Parameters and Results}
We now investigate the minimum number of parameters and their values required to display the original behaviour.  Firstly, the parameters for alive and injured friendly and enemy ISAACAs are made equal such that there is no distinction in behaviour of an alive or injured ISAACA.  Personality weights towards friendly forces (and thus clustering effects) were set to zero and found not to affect the overall ``front forming'' behaviour.  Switching the Advance parameter off also had no effect.  The final parameter values used are in Table \ref{tab:ISAACParametersAll}.
\begin{figure}[H]
  \centering
  \includegraphics[width=8.6cm]{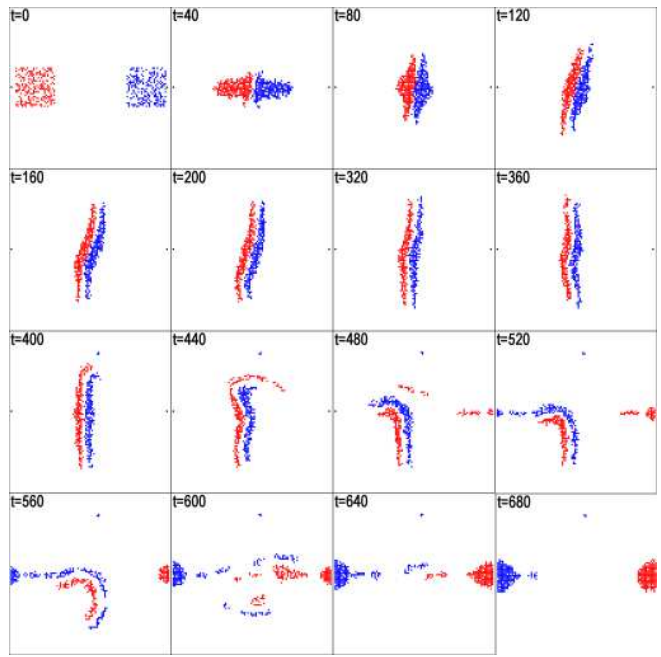}
  \caption[ISAAC Classic Fronts Scenario Screenshots]%
  {These snapshots were generated using the provided input file and demonstrates equivalent behaviour to those in \cite{andy97}.}
  \label{fig:ISAACFluid}
\end{figure}
\subsubsection{PDE Modification}
Our original equations are modified to reflect the addition of the combat meta-personality.  This is effectively a modification of the velocity term $\mathbf{f_{vel}}$ as we wish for only the direction of movement to be affected by the presence of the opposite force.  Diffusion, inter-force attraction and repulsion are to remain unaffected.
\begin{equation}
	\mathbf{f_{vel}} = \nabla \cdot\{u(\mathbf{C_u}u + A_{a}(K_{a} \ast u) - A_{r}u(K_{r} \ast u))\}
\end{equation}
We propose the following form of the velocity term, taking force $u$ as an example:
\begin{equation}
	\mathbf{C} = \left\{ 
		\begin{array}{l l}
			  C & \ if\ (\iint_{R}u\;dxdy-\iint_{R}v\;dxdy) > \Delta c\\
			  \\
			  -C & \ if\ (\iint_{R}u\;dxdy-\iint_{R}v\;dxdy) \leq \Delta c\\
		\end{array} \right. \label{eq:CFNewAdvection}
\end{equation}
where $(\iint_{R}u dxdy)$ is the number of friendly forces and $(\iint_{R}v dxdy)$ the number of enemy forces within the circular domain with radius corresponding to the sensor range $r_S$.  $C$ is a constant and represents the overall movement toward the $u$ or friendly flag.  If a numerical advantage greater than $\Delta c$ is reached, the force will proceed with a constant velocity toward its own flag.  Should this numerical advantage not be attained, the velocity of the force is reversed such that it will retreat away from its flag.
\subsubsection{PDE Results}
\begin{figure}[H]
	\centering
	\includegraphics[width=8.6cm]{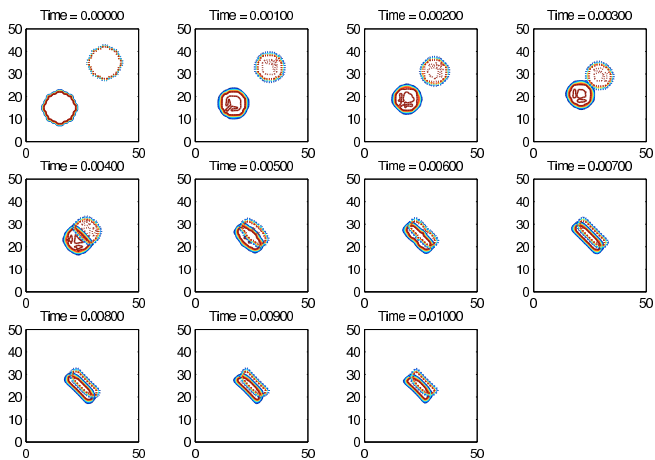}
	\caption{Our model: Classic Fronts.\label{fig:EinFluidAll1}}
\end{figure}
ISAAC has two areas of inbuilt randomness - new position selection when there are multiple positions with equal penalty function values, and the calculations of casualties/fatalities.  This randomness results in slight differences in the distribution of the two forces whereas the non-random continuous version has equal distributions.  The ISAAC results in Figure~\ref{fig:ISAACFluid} shows the forces ``slipping'' around each other due to these slight variations in distributions and then proceeding to their respective goals.  Figure~\ref{fig:EinFluidAll1} shows that our continuous version also forms long thin fronts that remain stationary.  Due to the steady reduction in both forces due to firing effects, their densities gradually decline however the position of the fronts remain stationary.  These tactics represent a classic style of attrition warfare.  If the simulation is allowed to progress, both forces will eventually decline to zero.  There is no randomness present to generate the slight differences or asymmetry in spatial distribution that leads to the forces manoeuvring around each other as in Figure~\ref{fig:ISAACFluid}.

	In order to introduce an asymmetry to the continuous model so as to mimic the asymmetry seen in the ISAAC results, the initial positions of each force distribution is changed such that the velocity vectors are offset.  That is, the forces will no longer collide ``head on''.
\begin{figure}[H]
	\centering
	\includegraphics[width=8.6cm]{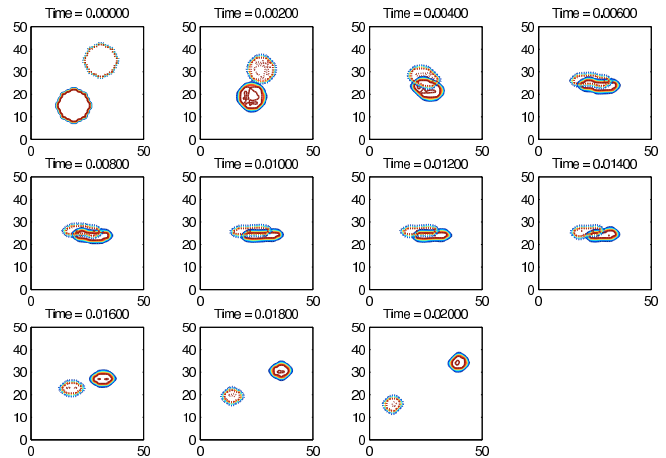}
	\caption{Our Model: Classic Fronts comparison with forces initially offset.\label{fig:EinFluidAll2}}
\end{figure}
\begin{figure}[H]
	\centering
	\includegraphics[width=6cm]{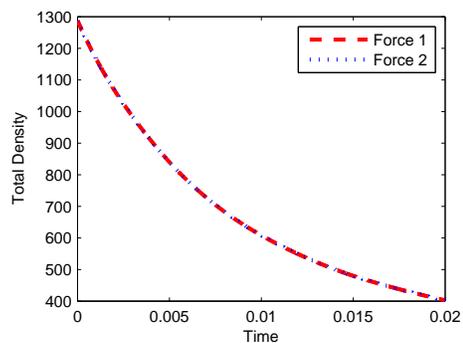}
	\caption{Losses for Figure~\ref{fig:EinFluidAll2}\label{fig:EinFluidOffsetLosses}}
\end{figure}
Again the formation of fronts is present however they form at a slight angle rather than vertically due to the initial offset in position.  A slight oscillation in the movement of each force is visible in Figure~\ref{fig:EinFluidAll2} as they slowly rotate around each other.  This is more noticeable in the movie of this scenario.  Once the forces pass, circular profiles with the respective desired minimum densities are reformed.  Comparing Figure~\ref{fig:ISAACFluid} to Figure~\ref{fig:EinFluidAll2}, the tight troop formation of the ISAAC agents seen prior to the forces making contact in Figure~\ref{fig:ISAACFluid} is not seen after the confrontation.  In the ISAAC results, troops stream to their respective flag in single file rather than proceeding in formation due to the effect of increasing flag proximity ($w_6$) on the penalty function.  However, maintaining a coherent profile throughout the \emph{entire} simulation, especially when a force is not in contact with another, is a highly desirable feature as a long thin profile may have an increased vulnerability to attack.
\subsection{Precess Scenario}
\subsubsection{Description}
This scenario is described by Ilachinski as showing a \textit{simple example of an emergent behaviour}.  Red and blue forces are quite different in personality, with red preferring to remain in close proximity while blue actively seeks the red force within sensor range.  Initially the rapidly advancing blue force is loosely formed and the slower red force advances in a higher density and thus smaller formation.  As the forces enter into each other's sensor range, the blue force partially surrounds the red and a slow precession begins.  The red force slowly continues advancing toward its flag, constantly pursued by the blue force as shown in Figure~\ref{fig:ISAACPrecessOffset}.  This precession behaviour, or rotation of both forces around an axis, as shown in the published ISAAC results and Figure~\ref{fig:ISAACPrecessOffset} occurs in only a small percentage of the test runs.  When precession was observed, both clockwise and anti-clockwise precession occurred with the forces often colliding with the domain boundaries.  Offsetting the position of the blue flag slightly to the left or right did result in a higher frequency of precession observed, with precession directionality depended upon the flag offset.

	Without the flag offset the usual behaviour is as follows using the initial conditions as given in Table \ref{tab:ISAACParametersAll}.  As the forces come into sensor range, the blue force surrounds the red until the majority of the blue force is located on the south western side of the red force.  As the personality weights for engaging in combat is stronger than for continuing to the blue flag, the blue force actively pursues the tightly clustered red force as it moves toward its red flag goal.  Both forces remain in coherent formations throughout the duration of the scenario as seen in Figure~\ref{fig:ISAACPrecess}.
	
	We propose that, similar to the Classic Fronts scenario, this precession behaviour arises due to the spatial asymmetry of the forces.  Depending on the spatial distribution on falling into sensor range, the precession will be either clockwise or anticlockwise, and offsetting the initial positions exacerbates this asymmetry thus increasing the frequency of precession occurring.
\subsubsection{ISAAC Parameters and Results}
Unlike the Classic Fronts scenario, no input file was provided by the developer of ISAAC with the software installation, only a series of snapshots with a partially complete list of parameter values.  Those unspecified values were determined through a trial and error procedure until the behaviour shown in \cite{andy97} could be reproduced.  In a similar manner as for the Classic Fronts scenario, we set the \emph{Advance} and \emph{Minimum Distance} parameters to zero and find that precession behaviour is still produced although more infrequently, suggesting that these parameters are unnecessary for this behaviour to occur.
\begin{figure}[H]
  \centering
  \includegraphics[width=8.6cm]{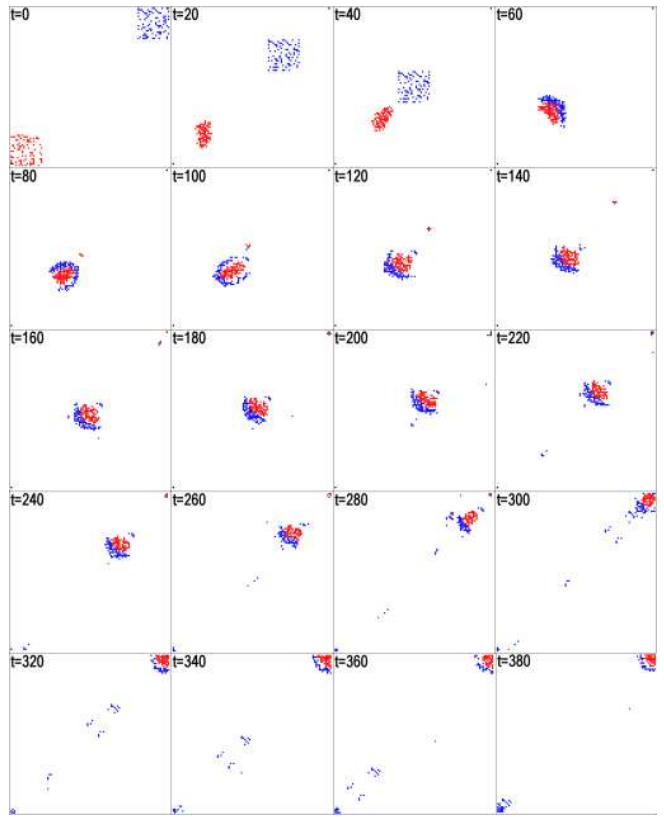}
  \caption{ISAAC Precess Scenario Screenshots}
  \label{fig:ISAACPrecess}
\end{figure}
Figure~\ref{fig:ISAACPrecess} shows the usual results found when running the input file as detailed in Table~\ref{tab:ISAACParametersAll}.  Note the absence of precession.  As indicated above, precession in both directions can be observed with the direction determined by the spatial asymmetry generated in that particular instance.
\begin{figure}[H]
	\centering
	\includegraphics[width=6cm]{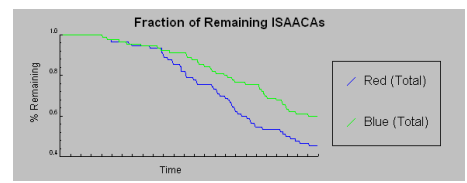}
	\caption{Losses for Figure~\ref{fig:ISAACPrecess}\label{fig:ISAACPrecessLosses}}
\end{figure}
\begin{figure}[H]
  \centering
  \includegraphics[width=8.6cm]{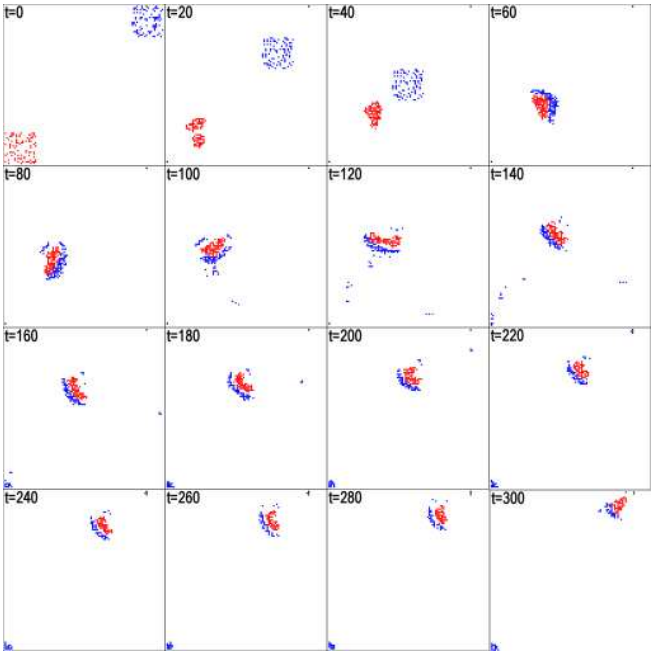}
  \caption{ISAAC Precess Scenario Screenshots with artificial Offset}
  \label{fig:ISAACPrecessOffset}
\end{figure}
By introducing an artificial offset, precession is almost always observed (Figure~\ref{fig:ISAACPrecessOffset}) and is dictated by the direction of the given offset.
\subsubsection{PDE Modifications}
We now seek to modify the velocity term in a similar way as with the Classic Fronts scenario.  Our modification to the velocity term described in (\ref{eq:CFNewAdvection}) does not adequately take into account the differences in offensive/defensive personalities as described in the ISAAC literature and must be expanded in order to include this.  Firstly, we define the number of friendly ($N_u$) and enemy forces ($N_v$) within the circular domain of radius $r_S$ (\emph{Sensor} range), used to determine the \emph{Combat} constraint by:
\begin{equation}\label{eq:ISAACFriendlyEnemyNumbers}
		N_v = \iint_{R}v\;dxdy;\ \ \ N_u = \iint_{R}u\;dxdy;
\end{equation} 
In order to allow for the two different types of attacking personalities, a switch $attack$ in the form of an integer of value $1$ or $-1$ is included.  This allows for a distinction between two force types - an aggressive and a defensive force.  An aggressive force will move toward the enemy regardless of superiority in numbers, while a defensive force will advance only with superior numbers and will otherwise retreat.  We now alter the velocity term again and arrive at the following form:
\begin{equation}\label{eq:PrecessAdvection}
	\mathbf{C} = \left\{ 
	\begin{array}{l l}
		C + \iint_{R}v\;dxdy & \ \mbox{$N_u-N_v \geq \Delta c$}\\ \\
		C + attack \times \iint_{R}v\;dxdy & \ \mbox{$N_u-N_v < \Delta c$}\\
	\end{array} \right.
\end{equation}
With a numerical superiority above the given $\Delta c$, the force has an additional attraction up the gradient of greatest concentration of the opposing force.  If this superiority level is not reached, this attraction is reversed, becoming a repulsion down the gradient of greatest concentration.  This effectively mimics the Combat constraint of ISAAC as described in Section \ref{sec:ISAACCombatConst}.
\subsubsection{PDE Results}
\begin{figure}[H]
	\centering
	\includegraphics[width=8.6cm]{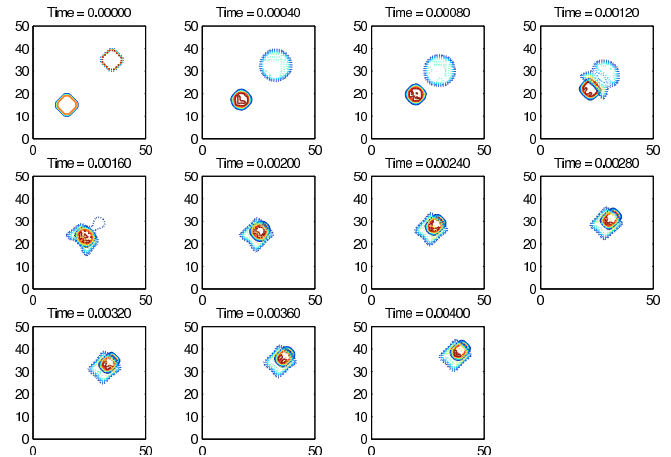}
	\caption{Our Model: Precess approximation.\label{fig:Precess}}
\end{figure}

Figure~\ref{fig:Precess} shows the first comparison at the ISAAC precession results and is markedly similar to Figure~\ref{fig:ISAACPrecess}.  Upon falling into the \emph{Sensor} range of the smaller footprint higher density Red force, the larger footprint lower density Blue force rapidly moves to surround it as $attack = 1$ for Blue.  As with the ISAAC scenario, this aggressive force pursues Red for the remainder of the simulation with the majority of the force located on the south western side of the Red force and a small portion of Blue located on the north eastern side of the Red force.
\begin{figure}[H]
	\centering
	\includegraphics[width=6cm]{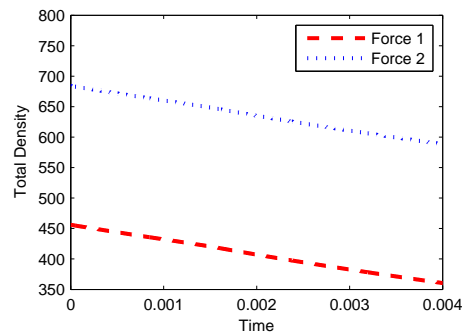}
	\caption{Losses for Figure~\ref{fig:Precess}\label{fig:PrecessLosses}}
\end{figure}
Note that no precession behaviour is seen which corresponds to the majority of ISAAC Precess scenario simulation results.  We now offset the initial positions of the forces as we did with the ISAAC scenario such that they pass to the left or right of each other to ascertain whether both clockwise and anticlockwise precession behaviour will be achieved in our continuous model.  We begin by attempting to induce anticlockwise precession.
\begin{figure}[H]
	\centering
	\includegraphics[width=8.6cm]{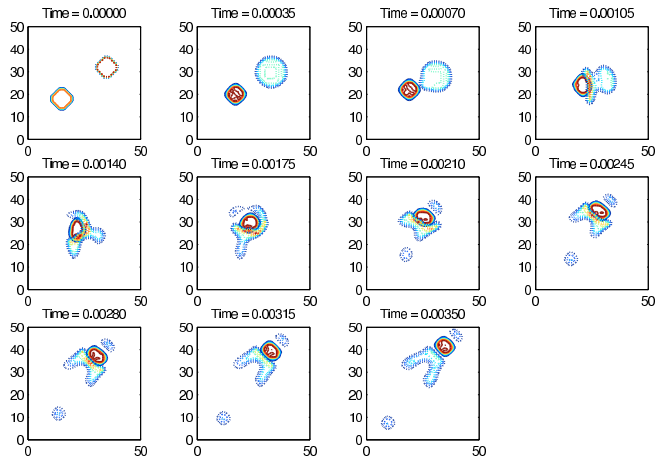}
	\caption{Our Model: Precess comparison with initial positions offset. Anticlockwise precession.  \label{fig:PrecessOffsetAnti}}
\end{figure}
\begin{figure}[H]
	\centering
	\includegraphics[width=6cm]{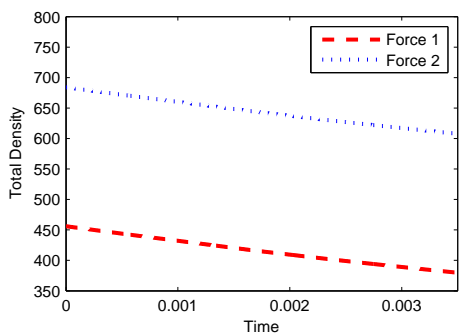}
	\caption{Losses for Figure~\ref{fig:PrecessOffsetAnti}\label{fig:PrecessOffsetAntiLosses}}
\end{figure}
	As expected, precession occurs due to the spatial asymmetry introduced into the scenario through the offset of the initial spatial distributions.
	The behaviour is markedly similar to that seen in Figure~\ref{fig:ISAACPrecessOffset}.  Upon falling into sensor range, the Red force is partially surrounded by Blue and an anti-clockwise precession occurs.  The Red force is then closely pursued by the majority of the Blue force following and a small portion directly in the path of Red.  Setting the initial positions with the reverse offset results in a clockwise precession. 
	
	We propose that, similar to the Classic Fronts scenario, the precession behaviour seen in the ISAAC results (Figure~\ref{fig:ISAACPrecessOffset}) arises due to the spatial asymmetry of the forces.  Depending on the spatial distribution on falling into sensor range, the precession will be either clockwise or anticlockwise, and offsetting the initial positions of the forces exacerbates this asymmetry thus increasing the frequency of precession occurring.  An artificial initial spatial asymmetry in the continuous version is necessary before precession is observed.
	
	As a result of the modification of the advection term through the addition of the $attack$ variable, minimal overlap between the forces is observed which is in keeping with the ISAAC behaviour.  However the positive value of $attack$ for the Blue force is sufficient to give rise to pursuing of the Red force, yet the threshold value prevents any significant overlap.  It is the effect of the repulsion down the gradient of highest concentration of enemy density which prevents any significant overlap.
\subsection{Circle Scenario}
\subsubsection{Description}
\begin{figure}[H]
  \centering
  \includegraphics[width=8.6cm]{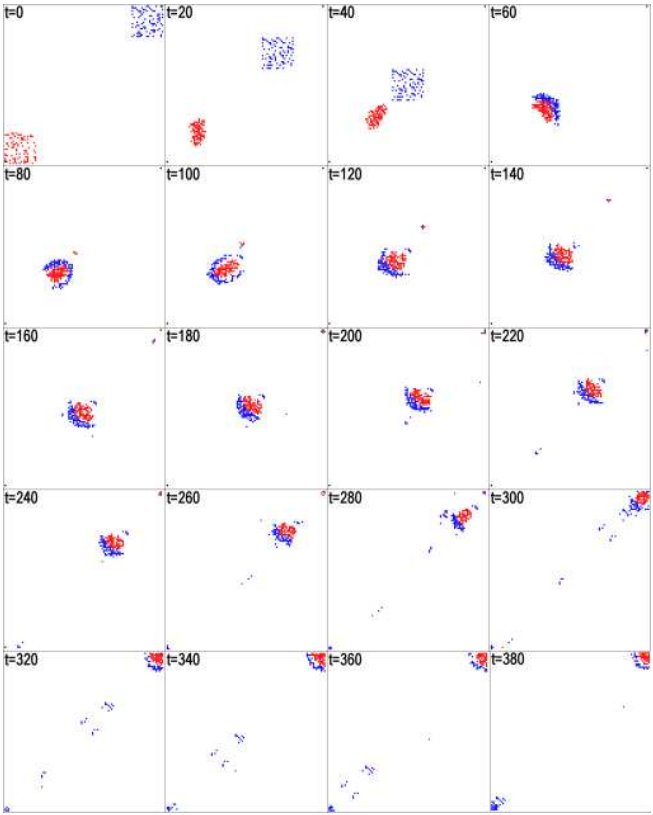}
  \caption{ISAAC Circle Scenario Screenshots\label{fig:ISAACCircle}}
\end{figure}

This scenario shows strikingly similar behaviour to the Precess scenario.  Both forces form dense formations with Blue slightly denser due to the higher cluster variable.  The more aggressive Red force envelops the Blue force entirely upon coming into sensor range and it is at this point, $time=140$, that the published snapshots stop.  This can give the impression that this formation remains static or stable after this time.  When running the simulation past this time step, the similarity with the Precess scenario from $time=120$ onwards becomes apparent.  Blue forces continue to move towards the Blue flag and Red forces begin to concentrate on the opposite side (Figure \ref{fig:ISAACCircle}, $time=200$).  For the remainder of the simulation, Red forces pursue the Blue force toward the Blue flag.  In this particular time series, a smaller section of the Blue force is separated and pursued by the Red force ($time=220$) in the same manner as the larger section.
\subsubsection{ISAAC Parameters and Results}
Unlike the Classic Fronts scenario and similar to the Precession scenario, no input file was provided with the software installation, only a series of snapshots with a partial list of parameters and corresponding values.  Initial distribution position and size values were not provided.  Numerous test runs using the given parameters showed a less dense blue force advancing more rapidly than the higher density red force.  Setting the \emph{Advance} and \emph{Minimum Distance} parameters to zero still produced the same behaviour.
\subsubsection{PDE Results}
Using similar parameter values to the Precess scenarios and using the same form of the advection term (\ref{eq:PrecessAdvection}), the following results are obtained.  The only differencea between the two scenarios are the fire parameters.
\begin{figure}[H]
	\centering
	\includegraphics[width=8.6cm]{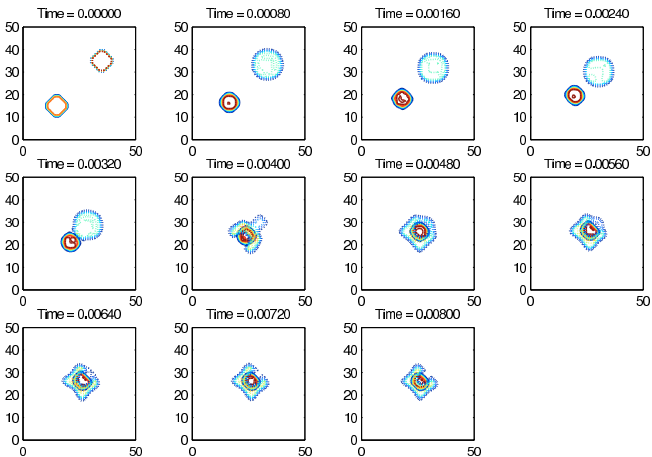}
	\caption{Our Model: Circle approximation.\label{fig:Circle1}}
\end{figure}
\begin{figure}[H]
	\centering
	\includegraphics[width=8.6cm]{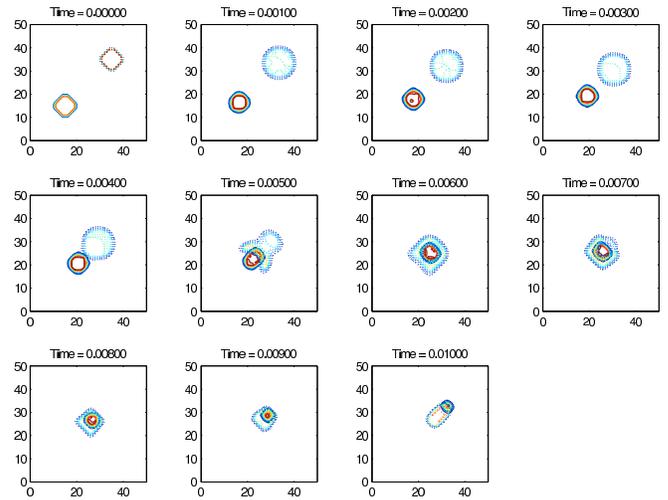}
	\caption{Our Model: Circle approximation with differing fire parameters.\label{fig:Circle2}}
\end{figure}
\begin{figure}[H]
  \centering
  \includegraphics[width=6cm]{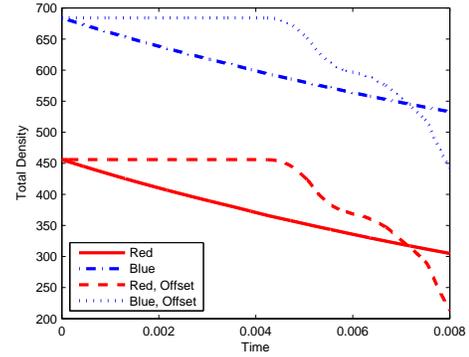}
  \caption{Losses for Figure~\ref{fig:Circle1} and Figure~\ref{fig:Circle2}.\label{fig:CircleLosses}}
\end{figure}
Similarly to the ISAAC results shown in Figure~\ref{fig:ISAACCircle}, in Figure~\ref{fig:Circle1} the denser Red force is surrounded by the more aggressive Blue force.  While the repulsion due to force inferiority maintains a degree of separation between the two forces, some overlap is present at the edges.  Once the Red force is surrounded, overall movement slows to almost a halt.  Density is gradually lost due to attrition and there is a slight imbalance in the Blue force's distribution around Red.  The slow advance of the Red force and imbalance in the Blue force distribution results in an eventual breakaway of Red.  The higher attrition rate used in Figure~\ref{fig:Circle2} shows that this almost stationary phase has a much shorter duration as attrition speeds the generation of an imbalance in the surrounding force's density.

	Note the nonlinear form of density loss shown in Figure~\ref{fig:CircleLosses}.  Lauren states that casualty rates are uneven on a turbulent nonlinear battlefield.  This is held as an important distinction between complex adaptive and convectional combat models.  Here we have shown that a conventional model is indeed capable of producing intermittent density losses.
\section{Discussion}
Cellular automata can be difficult to use for understanding the underlying dynamics of combat as stochasticity can hinder the extraction of conclusions from a model.  All scenarios presented here highlight the dangers associated with attributing \emph{intelligent} reasoning to behaviour shown, when this can be explained quite simply through the effects of the terms in our equations (\ref{eq:NewEquationForU}) and (\ref{eq:NewEquationForV}) and the spatial distribution of forces.  This can be seen quite simply in the Classic Fronts scenario.  Ilachinski describes the occurrence of the forces passing by one another as the agents having \emph{``found'' a way to sneak around}.  We believe it is the differences in distribution caused the inherent randomness of ISAAC that results in agents located at the end of the distributions becoming sufficiently distanced so that repulsive forces from enemy agents greatly reduces with respect to the attractive forces of the goal.  This dominant attractive force then results in those agents progressing towards the goal.  Although the attraction/repulsion of goals in our equations is not weighted with respect to separation distance as it is in the ISAAC penalty function, a slight offset or asymmetry of initial distributions demonstrates this same type of behaviour.

	This can also be seen clearly in the Precess scenario where the infrequently observed precession behaviour arises from asymmetries of the force profiles upon commencing combat, rather than being an unexplainable emergent behaviour.  The majority of ISAAC Precession simulations show behaviour similar to those obtained from the continuous equivalent without the initial distribution offset as shown in Figure \ref{fig:ISAACPrecess}.  Again the observed asymmetries are due to the stochasticity of ISAAC present in the movement and attrition algorithms, and are responsible for the precession seen in Figure \ref{fig:ISAACPrecessOffset} and \cite{andy97}.  When mimicked in our continuous equations through initially offset distributions, precession was observed.  In this case the underlying behaviour without the effects of stochasticity highlighted the similarity between the Precession and Circle scenarios.  Again Ilachinski infers a degree of intelligence or emergence by stating that there are \emph{a few stray `leakers' and an occasional group of a few Blue ISAACAs that choose[s] to leave the main battle and head toward Red's flag}.  This is similar to the comments made for the Classic Fronts scenario and our explanation for that behaviour - differences in distribution result in agents becoming sufficiently distanced from enemy agents such that the goal terms in the penalty function become significant.
	
	Our deterministic approach encapsulates basic motivational factors and demonstrates a variety of spatial behaviours.  Cohesive troop movement has only been achieved artificially in previous work: by using a desired initial distribution, low diffusion constant, and a sufficiently short overall simulation time, excessive diffusion is prevented and the troop profile cannot diffuse to a unrealistic spread.  We have demonstrated that by using relatively simple and physically meaningful form of partial differential equations, cohesive troop movement can be achieved and maintained, even when suffering loss of density through fire.  Forces and firing coefficients remain homogeneous, a criticism of the traditional Lanchester approach, yet the nonlinear nature of the equations are able to mimic those seen in ISAAC that have been labelled as complex.  A continuum of forces is able to behave in a manner similar to a collection of individual autonomous agents, and shows decentralised self-organisation and adaptation of tactics to suit a variety of combat situations.  This is a significant step toward developing a set of realistic continuous equations for combat modelling.
	
	Lauren \cite{l99} states that complex adaptive models of combat, such as ISAAC or MANA, display a rich variety of behaviour, a battlefield that is no longer linear and agent adaptivity.  How exactly do the agents perform this adapting?  The evolution of each agent's position is determined by a penalty function that remains unchanged throughout the entire simulation.  There is a danger in the anthropomorphisation of agents, insinuating agents have a reasoning and planning ability when this is obviously not the case.  Also these types of wargames concentrate heavily on the addition of extra communication ability between agents, shifting the emphasis to global or increasingly complex nonlocal features.  MANA includes many more states and subsequently many more triggers need to be defined to facilitate switching between these states.  Increasing the number of required parameters can cloud the process of deriving insight, a danger which has been shown in many of the basic scenarios presented here.  For example the removal of parameters such as the \emph{Advance} constraint in the Precess scenario did not prevent the reproduction of similar precession behaviour.  Lauren also states that \emph{conventional combat models behave largely as a series of attrition-driven fights} and \emph{ISAAC entities will only fight if conditions are suitable}.  If our model is viewed as a conventional combat model as it is based on using Lanchester firing terms, it can be argued that it also behaves as an ISAAC model due to the form of the Spatial Dynamics terms ($\mathbf{f_{diff}}$ and $\mathbf{f_{vel}}$).  Conversely it could be argued that ISAAC is a type of conventional combat model as the movement of each agent is determined in a comparable way to our continuous model with the inclusion of randomness.
	
	Considering the popularity of cellular automata based wargames in military and complex adaptive systems research, it is imperative that complementary avenues of research are undertaken in order to gain a greater insight into the nature of combat.  For a research topic such as combat modelling, emphasis of rare events or a misunderstanding of observed behaviour can have significant consequences.  This research has demonstrated that seemingly emergent and intelligent behaviour of agents sometimes seen only in atypical scenarios, can be explained through the numerical analysis of their continuous equivalent.  The need for many multiple agent simulations and the application of data mining techniques becomes much reduced with the simultaneous use of a continuous model.  MANA, much like ISAAC, can require approximately 600 runs to establish a mean result for some scenarios \cite{lau02}.  Using the approach employed here it is much easier to establish a mean behaviour with our continuous form.
	
	By treating weapons systems in the simple forms as presented here, behaviour is less likely to be obscured by potentially highly nonlinear interaction terms or weapons effects.  This could also lead to the exploration of the inherent nonlinearities in the Spatial Dynamics terms, which essentially describes the movement or tactics employed by a force.  For example, an increase in weapon lethality (kill probability) may not necessarily have the expected corresponding increase in casualties due to the effects of the enemy's tactics.  It may not be sound to assume a doubling of this probability will yield double the number of casualties.
	
	We suggest that our continuous model be used in conjunction with agent-based wargames to act as a combined testbed for the purpose of concept exploration.
\section{Future Work}
This research did not address other ISAAC scenarios such as LOCALCMD or GLBALCMD where communication between agents within a user-defined range affects the evolution of the scenarios.  Inclusion of the local command and global command functionalities and of command personality in the penalty function (\ref{eq:ISAACPenaltyFunction}) would require further significant modification to our continuous model.  Addition of these comparisons will result in a more complete development of a continuous counterpart to ISAAC.

	We have expanded our model to include a density response tactic where the equilibrium interior density is dependent upon the enemy density within the given sensor range.  This is similar to predator avoidance modelling in CA fish schooling modelling in \cite{va97} and has shown interesting dynamical changes in the test case scenarios thus far.  We envisage publishing these results in the near future.
\begin{acknowledgments}
I would like to acknowledge my supervisors Prof. Jim Franklin and Assoc. Prof. Gary Froyland for their guidance and assistance in preparing this manuscript.
\end{acknowledgments}

\appendix
\section{Parameters Used for Simulations}
\begin{table*}
\caption{\label{tab:ISAACParametersAll}ISAAC Parameters for Scenario Screenshots.  All other parameters are set to zero or \emph{No}.}
\begin{ruledtabular}
\begin{tabular}{ccccccccc}
	&\multicolumn{2}{c}{Classic Fronts (CF)}&\multicolumn{2}{c}{Precess (P)}&\multicolumn{2}{c}{Precess Offset (PO)}&\multicolumn{2}{c}{Circle (C)}\\
  Parameter&Red&Blue&Red&Blue&Red&Blue&Red&Blue\\ \hline
Squad Size&225&225&90&90&90&90&200&200 \\ 
$w_1$&0&0&25&10&25&10&10&25 \\
$w_2$&50&50&10&35&10&35&50&25 \\
$w_3$&0&0&75&10&75&10&0&75 \\
$w_4$&50&50&25&80&25&80&100&25 \\
$w_5$&0&0&0&0&0&0&0&0 \\
$w_6$&5&5&50&50&50&50&25&75 \\
$r_S$&5&5&5&5&5&5&5&5 \\
$r_F$&3&3&3&3&3&3&3&3 \\
$r_T$&2&2&2&3&2&3&3&3 \\
$w_M$&1&1&1&1&1&1&1&1 \\
Prob Hit&$2\times10^{-3}$&$2\times10^{-3}$&$2\times10^{-3}$&$2\times10^{-3}$&$2\times10^{-3}$&$2\times10^{-3}$&$10^{-3}$&$10^{-3}$\\
Max Sim tgts&5&5&All&All&All&All&999&999\\
Defence Measure&1&1&1&1&1&1&1&1\\
Cluster&0&0&10&3&10&3&3&15\\
Advance&0&0&0&0&0&0&0&0\\
Combat&3&3&4&-5&4&-5&-7&5\\
Initial Dist Centre x&15&85&10&90&10&90&10&90\\
Initial Dist Centre y&50&50&10&90&10&90&10&90\\
Size x&25&25&20&20&20&20&20&20\\
Size y&25&25&20&20&20&20&20&20\\
Flag x&1&99&1&99&1&90&1&99\\
Flag y&50&50&1&99&1&99&1&99 \
\end{tabular}
\end{ruledtabular}
\end{table*}

\begin{table*}
\caption{\label{tab:PDEParametersAll}PDE Parameters.}
\begin{ruledtabular}
\begin{tabular}{ccccccccccccc}
	&\multicolumn{2}{c}{CF}&\multicolumn{2}{c}{CFO}&\multicolumn{2}{c}{P}&\multicolumn{2}{c}{Anticlock P}&\multicolumn{2}{c}{C}&\multicolumn{2}{c}{C, High $d$}\\
  Parameter&$u$&$v$&$u$&$v$&$u$&$v$&$u$&$v$\\ \hline
$\mathrm{ID}$&8&8&8&8&8&12&8&12&8&12&8&12\\ 
$\rho$&5&5&5&5&5&5&5&5&5&5&5&5\\
$\mu$&(15,15)&(35,35)&(19,15)&(31,35)&(15,15)&(35,35)&(15,18)&(35,32)&(15,15)&(35,35)&(15,15)&(35,35)\\
$\mathrm{IT}$&0.5&0.5&0.5&0.5&1&1&1&1&1&1&1&1\\
$r_{a,r}$&5&5&5&5&5&5&5&5&5&5&5&5\\
$D$&5&5&5&5&5&5&5&5&5&5&5&5\\
$C$&(20,20)&(-20,-20)&(20,20)&(-20,-20)&(60,60)&(-60,-60)&(60,60)&(-60,-60)&(20,20)&(-20,-20)&(20,20)&(-20,-20)\\
$A_a$&5&5&5&5&5&5&5&5&5&5&5&5\\
$A_r$&0.5&0.5&0.5&0.5&0.5&1&0.5&1&0.5&1&0.5&1\\
$r_S$&3&3&3&3&3&7&3&7&3&4&3&5\\
$\Delta c$&100&100&100&100&$10^6$&4&$10^6$&4&18&18&20&20\\
attack& & & & &-1&-1&-1&-1&-1&-1&-1&-1\\
$d$&$2\times10^{-6}$&$2\times10^{-6}$&$2\times10^{-6}$&$2\times10^{-6}$&$2\times10^{-6}$&$2\times10^{-6}$&$2\times10^{-6}$&$2\times10^{-6}$&$10^{-5}$&$10^{-5}$&$10^{-4}$&$10^{-4}$\\
$\beta$&$8\times10^{-8}$&$8\times10^{-8}$&$8\times10^{-8}$&$8\times10^{-8}$&$8\times10^{-8}$&$8\times10^{-8}$&$8\times10^{-8}$&$8\times10^{-8}$& & & & \\
$\nu$&0.2&0.2&0.2&0.2&0.2&0.2&0.2&0.2& & & & \\
$\tau(t=0)$&$10^{-7}$& &$10^{-7}$& &$10^{-7}$& &$10^{-7}$&&$10^{-7}$& &$10^{-7}$&\\
end time&$10^{-2}$& &$2\times10^{-2}$& &$4\times10^{-3}$& &$3.5\times10^{-3}$&&$8\times10^{-3}$& &$10^{-2}$&\\
atol&$10^{-3}$& &$10^{-3}$& &$10^{-3}$& &$10^{-3}$& &$10^{-3}$& &$10^{-3}$&\\
rtol&$10^{-3}$& &$10^{-3}$& &$10^{-3}$& &$10^{-3}$& &$10^{-3}$& &$10^{-3}$&\\
\end{tabular}
\end{ruledtabular}
\end{table*}

\end{document}